\documentclass[oneside]{article}
\usepackage{amsmath,amssymb,amscd,amsthm,verbatim,alltt,amsfonts,array}
\usepackage{mathrsfs}
\usepackage[english]{babel}
\usepackage{latexsym}
\usepackage{amssymb}
\usepackage{euscript}
\usepackage{graphicx}
\usepackage{arydshln}
\usepackage[dvipsnames]{xcolor}

\newtheorem{theorem}{Theorem}
\theoremstyle{plain}

\numberwithin{equation}{section}

\begin{document}

{\footnotesize%
\hfill
}

  \vskip 1.2 true cm

\begin{center} {\bf Homological monotonicity for configuration spaces of manifolds} \\
          {by}\\
{\sc Muhammad Yameen}
\end{center}

\pagestyle{myheadings}
\markboth{Homological monotonicity for configuration spaces of manifolds}{Muhammad Yameen}

\begin{abstract}
Consider the configuration spaces of manifolds. An influential theorem of McDuff, Segal and Church shows that the (co)homology of the unordered configuration space is independent of number of points in a range of degree called the stable range. We study the another important (and general) property of unordered configuration spaces of manifolds (not necessarily orientable, and not necessarily admitting non-vanishing vector field) that is \emph{homological monotonicity} in unstable part. We show that the homological dimension of unordered configuration spaces of manifolds in each degree is monotonically increasing. Our results show that the monotonicity property is not depend on the differential structure and orientiability of manifold. 
\end{abstract}

\begin{quotation}

\end{quotation}

\thispagestyle{empty}

\section{Introduction}

\label{sec:intro}


For any manifold $M$, let
$$F_{k}(M):=\{(x_{1},\ldots,x_{k})\in M^{k}| x_{i}\neq x_{j}\,for\,i\neq j\}$$
be the configuration space of $k$ distinct ordered points in $M$ with induced topology. The symmetric group $S_{k}$ acts on $F_{k}(M)$ by permuting the coordinates. The quotient $$C_{k}(M):=F_{k}(M)/S_{k} $$
is the unordered configuration space with quotient topology. All manifolds are assumed to be finite type (Betti numbers are finite numbers). The boundary of manifold $M$ is written $\partial M.$ Throughout the paper, $U$ is a closed subset of $M$ such that $\partial M\cap U=\phi,$ and $M-U$ is connected. 

An influential theorem of McDuff \cite{MD}, Segal \cite{S} and Church \cite{C} shows that the (co)homology of the unordered configuration space $C_{k}(M)$ is independent of number of points $k$ in a range of degree called the stable range:
\begin{theorem}\label{Mono1}
Let $M$ be a manifold. For each $i\geq0,$ the map
$$k\rightarrow \emph{dim} H_{i}(C_{k}(M);\mathbb{Q})$$ 
is eventually constant.
\end{theorem}
Ellenberg and Wiltshire-Gordon \cite{EW} proved that the homological dimension of unordered configuration spaces of orientable manifolds admitting non-vanishing vector field is non-decreasing. We prove that the conditions of orientablility and admitting vanishing vector field are not necessary.

\begin{theorem}\label{Mono1}
Let $M$ be a non-orientable manifold. For each $i\geq0,$ the map
$$k\rightarrow \emph{dim} H_{i}(C_{k}(M);\mathbb{Q})$$ 
is monotonically increasing.
\end{theorem}
\begin{theorem}\label{Mono2}
Let $M$ be a orientable manifold and $\partial M\cup U\neq\phi.$ For each $i\geq0,$ the map
$$k\rightarrow \emph{dim} H_{i}(C_{k}(M-U);\mathbb{Q})$$ 
is monotonically increasing.
\end{theorem}

\subsection*{General conventions}

$\bullet$ We work throughout with finite dimensional graded vector spaces. The degree of an element $v$ is written $deg(v)$.\\\\
$\bullet$ The symmetric algebra $Sym(V^{*})$ is the tensor product of a polynomial algebra and an exterior algebra:
$$ Sym(V^{*})=\bigoplus_{k\geq0}Sym^{k}(V^{*})=Poly(V^{even})\bigotimes Ext(V^{odd}), $$
where $Sym^{k}$ is generated by the monomials of length $k.$\\\\
$\bullet$ The $n$-th suspension of the graded vector space $V$ is the graded vector space $V[n]$ with
$V[n]_{i} = V_{i-n},$ and the element of $V[n]$ corresponding to $a\in V$ is denoted $s^{n}a;$ for example
$$ H^{*}(S^{d};\mathbb{Q})[n] =\begin{cases}
      \mathbb{Q}, & \text{if $*\in\{n,n+d \}$} \\
      0, & \mbox{otherwise}.\\
   \end{cases} $$ \\\\
$\bullet$ We write $H_{-*}(M;\mathbb{Q})$ for the graded vector space whose degree $-i$ part is
the $i$-th homology group of $M;$ for example
$$ H^{-*}(S^{d};\mathbb{Q}) =\begin{cases}
      \mathbb{Q}, & \text{if $*\in\{-d,0.\}$} \\
      0, & \mbox{otherwise}.\\
   \end{cases} $$


\section{Chevalley–Eilenberg complex}
F\'{e}lix--Thomas \cite{F-Th} (see also \cite{F-Ta}) constructed a Sullivan model for the rational cohomology of configuration spaces of closed orientable even dimensional manifolds. Furthermore, Knudsen \cite{Kn} extended the result of F\'{e}lix--Thomas for general even dimensional manifolds. \\\\
Let us introduced some notations. Consider two graded vector spaces $$V^{*}=H^{-*}_{c}(M;\mathbb{Q}^{w})[d],\quad W^{*}=H_{c}^{-*}(M;\mathbb{Q})[2d-1]:$$
where
$$ V^{*}=\bigoplus_{i=0}^{d}V^{i},\quad W^{*}=\bigoplus_{j=d-1}^{2d-1}W^{j}.$$
We choose bases in $V^{i}$ and $W^{j}$ as 
$$V^{i}=\mathbb{Q}\langle v_{i,1},v_{i,2},\ldots\rangle,\quad W^{j}=\mathbb{Q}\langle w_{j,1},w_{j,2},\ldots\rangle$$
(the degree of an element is marked by the first lower index; $x_{i}^{l}$ stands for the product $x_{i}\wedge x_{i}\wedge\ldots\wedge x_{i}$ of $l$-factors). Always we take $V^{0}=\mathbb{Q}\langle v_{0}\rangle$. Now consider the graded algebra
$$ \Omega^{*,*}_{k}(M)=\bigoplus_{i\geq 0}\bigoplus_{\omega=0}^{\left\lfloor\frac{k}{2}\right\rfloor}
\Omega^{i,\omega}_{k}(M)=\bigoplus_{\omega=0}^{\left\lfloor\frac{k}{2}\right\rfloor}\,(Sym^{k-2\omega}(V^{*})\otimes Sym^{\omega}(W^{*})) $$
where the total degree $i$ is given by the grading of $V^{*}$ and $W^{*}$. We called $\omega$ is a weight grading. The differential $\partial:Sym^{2}(V^{*})\rightarrow W^{*}$ is defined as a coderivation by the equation 
$$\partial(s^{d}a\wedge s^{d}b)=(-1)^{(d-1)|b|}s^{2d-1}(a\cup b),$$ where $$\cup\,:H^{-*}_{c}(M;\mathbb{Q}^{w})^{\otimes2}\rightarrow H^{-*}_{c}(M;\mathbb{Q})$$
(here $H^{-*}_{c}$ denotes compactly supported cohomology of $M$). The degree of $\partial$ is $-1.$ It can be easily seen that $s^{d}a,\,s^{d}b\in V^{*}$ and $s^{2d-1}(a\cup b)\in W^{*
}.$ The differential $\partial$ extends over  $\Omega^{*,*}_{k}(M)$ by co-Leibniz rule. By definition the elements in $V^{*}$ have length 1 and weight 0 and the elements in $W^{*}$ have length 2 and weight 1. By definition of differential, we have 
$$\partial:\Omega^{*,*}_{k}(M)\longrightarrow\Omega^{*-1,*+1}_{k}(M).$$

\begin{theorem}
	If $d$ is even, $H_{*}(C_{k}(M);\mathbb{Q})$  is isomorphic to the homology of the complex 
	$$ (\Omega^{*,*}_{k}(M),\partial).$$
\end{theorem}

\section{Proofs of main results}
In this section, we give the proofs of Theorems \ref{Mono1} and \ref{Mono2}.\\\\
\textit{Proof of Theorems \ref{Mono1} and \ref{Mono2}.}\\\\
\textbf{Case I.} Let $d$ is odd. We have an isomorphism \cite{BCT}:
$$H_{*}(C_{k}(M);\mathbb{Q})\cong Sym^{k}(H_{*}(M;\mathbb{Q})).$$ 
We can write the decomposition of algebras: 
$$Sym^{k}(V^{*})\cong (\langle v_{0}\rangle\otimes Sym^{k-1}(V^{*}))\oplus Sym^{k}(V^{*})/\langle v_{0}\rangle,$$
here we take $V^{*}\cong H_{*}(M;\mathbb{Q})$ and $v_{0}$ is a corresponding element of $H_{0}(M;\mathbb{Q}).$
Moreover, we have an isomorphism :
$$H_{*}(C_{k-1}(M);\mathbb{Q})\cong \langle v_{0}\rangle\otimes Sym^{k-1}(V^{*}).$$
As a consequence, we get 
$$H_{*}(C_{k}(M);\mathbb{Q})\cong H_{*}(C_{k-1}(M);\mathbb{Q})\oplus Sym^{k}(V^{*})/\langle v_{0}\rangle.$$
This implies that for each $i\geq0,$ the map
$$k\rightarrow \emph{dim} H_{i}(C_{k}(M-U);\mathbb{Q})$$ 
is monotonically increasing.\\\\
\textbf{Case II.} Let $d$ is even and $M$ is not closed orientable. We have decomposition of algebras:
$$\Omega^{*,*}_{k}(M)\cong (\langle v_{0}\rangle\otimes \Omega^{*,*}_{k-1}(M))\oplus \Omega^{*,*}_{k}(M)/\langle v_{0}\rangle.$$
The subspace $\Omega^{*,*}_{k}(M)/\langle v_{0}\rangle$ is $\partial-$invariant:
$$\partial(\Omega^{*,*}_{k}(M)/\langle v_{0}\rangle)\subseteq \Omega^{*,*}_{k}(M)/\langle v_{0}\rangle.$$
The degrees of elements of $V^{*}$ are less than $d.$ Also, the degree of elements of $W^{*}$ are greater than and equal to $d-1.$ The differential $\partial$ has degree $-1.$ Therefore, $\partial(\langle v_{0}\rangle\otimes V^{*})=0.$ 
This implies that the subspace $\langle v_{0}\rangle\otimes \Omega^{*,*}_{k-1}(M)$ is also $\partial-$invariant:
$$\partial(\langle v_{0}\rangle\otimes \Omega^{*,*}_{k-1}(M))\subseteq \langle v_{0}\rangle\otimes \Omega^{*,*}_{k-1}(M).$$
As a consequence, we get the decomposition of complexes:
$$(\Omega^{*,*}_{k}(M),\partial)\cong (\langle v_{0}\rangle\otimes (\Omega^{*,*}_{k-1}(M),\partial)\oplus (\Omega^{*,*}_{k}(M)/\langle v_{0}\rangle,\partial).$$
As a consequence, we have an isomorphism :
$$H_{*}(C_{k}(M);\mathbb{Q})\cong H_{*}(C_{k-1}(M);\mathbb{Q})\oplus H_{*}(\Omega^{*,*}_{k}(M)/\langle v_{0}\rangle,\partial).$$
This implies that for each $i\geq0,$ the map
$$k\rightarrow \emph{dim} H_{i}(C_{k}(M-U);\mathbb{Q})$$ 
is monotonically increasing.
\section{Final remark}
Napolitano \cite{N} computed the first integral homology group of 2-sphere:
$$H_{1}(C_{k}(S^{2});\mathbb{Z}=\mathbb{Z}/(2k-2)\mathbb{Z}.$$ The integral homological stability fails for 2-sphere. Similar behavior in case of closed orientable manifolds is also observe for monotonicity. If $M$ is an orientable manifold (admitting non-vanishing vector field) and $\partial M\cup U=\phi,$ then the monotonicity is also hold for $C_{k}(M)$ (see \cite{EW}). However, this is not true in general for all closed orientable manifolds. For example the map $$k\rightarrow \emph{dim} H_{2d}(C_{k}(S^{2d});\mathbb{Q})$$
is not monotonically increasing (see \cite{RW}). The situation is unclear in this case, and we ask the following.\\\\
\textbf{Question.} Let $M$ be an orientable closed manifold with $\chi(M)\neq0.$ Then under what condition the map 
$$k\rightarrow \emph{dim} H_{i}(C_{k}(M);\mathbb{Q})$$ 
is monotonically increasing for each $i\geq0?$ \\\\

\noindent\textbf{Acknowledgement}\textit{.} The author gratefully acknowledge the support from the ASSMS, GC university Lahore. This research is partially supported by Higher Education Commission of Pakistan.

\vskip 0,65 true cm




\vskip 0,65 true cm







\null\hfill  Abdus Salam School of Mathematical Sciences,\\
\null\hfill  GC University Lahore, Pakistan. \\
\null\hfill E-mail: {yameen99khan@gmail.com}

\end{document}